\documentclass[10pt]{article}

 \usepackage{color,graphicx}
 \usepackage{amsmath,amssymb,amsthm,amsfonts,epic}

 \newtheorem{theorem}{Theorem}
 \newtheorem{prop}{Proposition}

 \usepackage{bm}
 \usepackage{bbm}

\begin{document}

\title{A polling system whose stability region
  depends on a whole distribution of service times}

\author{Natalia Chernova \and Sergey Foss \and Bara Kim}

\date{}

\maketitle

\begin{abstract}
We present an example of a single-server polling system with two
queues and an adaptive service policy where the stability region
depends on the expected values of all the primitives and also on a
certain exponential moment of the service-time distribution in one
of the queues. The latter parameter can not be determined, in
general, in terms of any finite number of power moments. It follows
that 
the fluid approximation approach may not be an appropriate tool for
the stability study of this model.


\end{abstract}

\section{Introduction} \label{section1}

We deal with a class of cyclic polling systems which are
single-server systems where the server visits a finite number of
queues in a cyclic order and serves customers there. The stability
and performance analysis of polling models has been a hot topic in 
the last 15-20 years. For overviews on the current progress in the
studies of polling models, see, for example, Borst (1995); Boxma et
al. (2009); Wierman et al. (2007); Winands et al. (2009); Boon
(2011). Stability conditions in many polling models and, in
particular, in those with state-independent service policies may be
analysed via the fluid approximation approach (see Rybko and
Stolyar (1992) 
and Dai (1995)) 
that is based on the Functional Strong Law of Large Numbers and
involves only the first moments of driving random elements. However,
it is known (may be, not sufficiently broadly) that, in general, a
dependence of stability conditions on distributions (of stochastic
characteristics of polling models) may be more complex and may
involve into consideration not only power moments but whole
primitive distributions.

In this article, we present a simple and {\it natural} example of a
polling system with an adaptive service mechanism. More precisely,
we introduce a 2-station system with an exhaustive service policy at
one of them and a limited policy at the other, where the limited
policy increases the limit level if the server finds the other
station empty. We assume the input streams to be Poisson, the
switch-over times to be exponential, while the service time
distributions in both queues are assumed to be rather general.

We show that the stability conditions depend on parameters of the
Poisson processes and of the exponential random variables, on the
first moments of the service time distributions from both stations,
and also on an exponential moment (!) of the service time at the
station where the limited policy is in use.
In particular, this shows that, in general, a knowledge of any
finite number of power moments is not enough to determine the
stability region. This also shows that the direct probabilistic
methods are still of great importance in the stability analysis of
stochastic networks and that the fluid approximation approach is not
as universal as many researchers think.

This paper may be viewed as complementary to Chernova et al.
(2011), where the authors analysed another single-server cyclic
polling system with three queues and a similar adaptive rule, and 
obtained upper and lower bounds for the stability region using the
fluid approximation approach. They provided the simulation results
showing that stability conditions for systems with the same power
moments of primitive random variables (up to the third moment), but
with different distributions, may also differ.

 \section{The Model and Its Stability Conditions}

 We consider a polling system with 2 stations and a single
 server that consequently visits the stations and serve customers there.
 Customers arrive at station $k=1,2$ in a Poisson
 stream of intensity $\lambda_k >0$. Service times at station 
 $k$ form an
 i.i.d. sequence having distribution function $B_k$ with positive finite
 mean $b_k$.  It takes an exponential time with mean $\gamma$ for the server to travel either from station 1 to station 2 or  
 from station 2 to station 1.
 We assume all primitive random variables (interarrival, service, and travel times)
 to be mutually independent.
 Server follows the exhaustive policy at station 2 and an adaptive
 limited policy at station 1. In more detail,
 the server works in ``cycles'', and each cycle starts when the server arrives at station 2.

 In more detail, if the server finds at least one customer 
 at station 2 when he 
 arrives there,
 he starts to serve customers one-by-one (including new
 arrivals) until he empties the queue. Then he travels to station 1 (during an exponential time), serves one customer there (if there
 is any) or leaves it immediately (if the queue is empty), 
 and then travels back to station 2 (during 
 another exponential time). 
 This is a {\it standard}  cycle.

 If, upon his arrival to station 2, the server finds it empty, the 
 next cycle is {\it modified}:
 the server is allowed to serve $m$ extra customers at station 1. More precisely,
 he starts the cycle with his travel to station 1; then he provides
 a service 
 at station 1 until either he completes $1+m$ services or the station becomes empty 
 (here $m$ is a fixed non-negative integer). Then he travels back to
 station 2.

 We need the cycles to be finite and, moreover, to have a finite mean. For that, we need the
 mean time to empty the second queue to be finite. So we assume the following condition to hold:
 \begin{equation}\label{stab1}
 \rho := \lambda_2 b_2 <1.
 \end{equation}

We are interested in {\it stability conditions} for this model. For 
$t\ge 0$, let $Q_k(t)$, $k=1,2$, be the queue length at
station $k$ at time $t$. Let $L(t)$ be the location of the server.
More precisely, $L(t)=k$ if the server is in service at station $k$,
and $L(t)=(1,2)$ or $L(t)=(2,1)$ if the server is travelling either
from $1$ to $2$ or back. Let $R(t)$ be the residual service time if
the server is in service (and let $R(t)=0$ if the server is travelling).
Let $I(t)$ be 1 if the cycle is standard and 2 if the cycle is
modified. Then the vector $Z(t)= (Q_1(t),Q_2(t),L(t),R(t), I(t))$
represents the current state of the system. We assume the process
$Z(t)$ to be right-continuous. This is a continuous-time Markov
process. We seek for conditions for the distribution of $Z(t)$ to
converge weakly, as $t\to\infty$, to a proper limiting distribution.

\subsection{Embedded Processes}

 For $i=1,2$ and $n=1,2,\ldots$, let $Q_{n,i}$ be the number of customers at station $i$ when
 the server arrives there for the $n$th time.
 Let $Z_n = (Q_{n,1},Q_{n,2})$, $n=1,2,\ldots$.
 Since the input processes are Poisson, the sequence $\{Z_n\}$
 forms a time-homogeneous Markov
 chain taking values in a countable state space ${\cal Z}_+^2$, and any  two states from the state space communicate to each other. 
 Since cycle lengths include exponential travel times, their distributions are absolutely continuous,
 and, in particular, the following conclusion may be deduced (see, e.g., Asmussen (2003) for the background):\\

\begin{prop}
 Under the condition \eqref{stab1},
 the process $\{Z(t):t\ge 0\}$ has a proper limiting distribution (is {\it stable}) if and only if the Markov Chain
 $\{Z_n:n=1,2,\ldots\}$ is {\it positive recurrent}.
 \end{prop}

 \subsection{Auxiliary Markov Chain}
 In this subsection, we assume that initially, at time $t=0$,
 there is infinitely many 
 customers at station 1 and no customers at station 2.
 Assume also that, at time $0$, the server starts to travel from 
 station 2 to station 1.
 For $t\ge 0$, let $\hat Z(t)$ be the number of customers at station 
 2 at
 time $t$ (so $\hat Z(0)=0$.
 Then $\{\hat Z(t):t\ge 0\}$ is a regenerative process, and a
 regeneration occurs 
 when the server arrives to queue 2 and finds
 it empty.
 The first renegerative time is denoted by $T>0$. 

 For $n=1,2,\ldots$, let $t_n$ be the $n$th return time of the server to station 2, so $t_0:=0$.
 Let  $\hat Q_n:= \hat Z(t_n)$, $n=0,1,2,\ldots$.
 Then $\{\hat Q_n:n=0,1,2,\ldots\}$ forms a Markov chain.
 Moreover,
 random sequence $X_{n} := {\bf I} (\hat Q_{n}>0)$ is also Markov.
 It takes only two values: $0$ if the queue is empty, and $1$,
 otherwise.
 Indeed, if $X_{n}=0$, then $X_{n+1}=0$  if there
 is no arrivals to station 2 during two travel times (from station 2 
 to station 1 and back) and $1+m$ service times
 at station 1, i.e.,
 $$
 p_{0,0} := {\mathbf P} (X_{n+1}=0 \ | \ X_n=0) = \left(\frac{\gamma^{-1}}{\lambda_{2}+\gamma^{-1}}\right)^{2}
 \left(B_{1}^{*}(\lambda_{2})\right)^{1+m}.
 $$
 Notice that $\frac{\gamma^{-1}}{\lambda_{2}+\gamma^{-1}}$ is the
 probability that there is no arrival to station 2 during a travel
 time of the server, and $B_{1}^{*}(\lambda_{2}) = \int_{0}^{\infty} e^{-\lambda_{2}t}
 dB_{1}(t)$ is the probability that there is no arrivals to station 2
 during a service time of a customer at station 1.
 Clearly $$B_{1}^{*}(\lambda_{2})={\mathbf E}
 e^{-\lambda_{2}\sigma_{1}},$$
 where $\sigma_1$ is a generic service time of a customer at
 station 1.
 For the Markov chain $\{X_n\}$, the transition probability from 0
 to 1 is given by
 $ p_{0,1}:= {\mathbf P} (X_{n+1}=1 \ | \ X_n=0) = 1-p_{0,0}$.
 If $X_n=1$, then, for any queue length, the server
 first empties queue 2. After that, the server travels to station
 1, serves a customer at station 1, and travels back to station 2.
 Therefore, given $X_n=1$, the probability that the server finds
 station 2 empty when returns there next time is
  $$
  p_{1,0}:= {\mathbf P} (X_{n+1}=0 \ | \ X_n=1) =
  \left(\frac{\gamma^{-1}}{\lambda_{2}+\gamma^{-1}}\right)^{2}
 B_{1}^{*}(\lambda_{2}).
 $$
 Finally, the transition probability $p_{1,1}$ is given by 
 $ p_{1,1}:= {\mathbf P} (X_{n+1}=1 \ | \ X_n=1) = 1-p_{1,0}$.

 Let $\nu$ be the first return time to 0 for the Markov chain
 $\{X_n\}$,  
 $$
 \nu := \min \{ n\ge 1 \ : \ X_{n}=0 \}.
 $$
 Then, by straightforward calculations, 
 $$
 {\mathbf E} \nu = 1 + \frac{p_{0,1}}{p_{1,0}}
 $$
 and
 $$ T = \sum_{i=1}^{\nu} \psi_{i},
 $$
where $\psi_{i}$ is a duration of the $i$th cycle (in continuous
time), which is the sum of the service times at both stations and of
 the travel times.
 Hence
\begin{eqnarray*}
  {\mathbf E} T &=& \sum_{i=1}^{\infty}
  {\mathbf E} (\psi_{i} {\mathbf I} (\nu \ge i)) \\
  &=&
  {\mathbf E} \psi_{1}+\sum_{i\ge 2}{\mathbf E}
  (\psi_{i}{\mathbf I} (\nu \ge i-1) {\mathbf I} (\hat Q_{i-1}\ge 1))\\
  &=&
  {\mathbf E} \psi_{1}+
  \sum_{i\ge 2}
  \left(
    {\mathbf E} (\hat Q_{i-1} | \nu \geq i-1)\frac{b_{2}}{1-\rho} {\mathbf P} (\nu \ge i-1)
    + (2\gamma + b_{1}) {\mathbf P}(\nu \ge i)\right),
  \end{eqnarray*}
  where we have used the fact that, in an $M/G/1$ queue with input rate $\lambda_2$
  and mean service time $b_2$, the mean of a busy period is 
  $b_{2}/(1-\rho )$.
  Note that
  $$
  {\mathbf E} \psi_{1}= 2\gamma + (1+m) b_{1},
  \quad
  {\mathbf E} \hat Q_{1}= (2\gamma  +
    (1+m)b_{1})\lambda_{2},
  \quad
  {\mathbf E} (\hat Q_{i-1} | \nu \geq i-1) =
  (2\gamma  +
    b_{1})\lambda_{2},
  $$
  for any $i\ge 3$.
By substituting these expressions, we obtain
  $$
  {\mathbf E} T = (2\gamma +b_{1})
  \frac{1}{1-\rho}{\mathbf E} \nu + \frac{1}{1-\rho} mb_{1}.
  $$
  Since the number of service completions at station 1 before time 
  $T$ is
  $\nu+m$, the mean service rate at station 1 is given by
  $$
  r:= \frac{{\mathbf E}[ \nu +m]}{{\mathbf E} T} =
  \frac{({\mathbf E} \nu +m)(1-\rho
    )}{(2\gamma +b_{1})
    {\mathbf E} \nu + mb_{1}}.
$$

  \subsection{Stability Analysis}

  Based on the calculations from the previous subsection, we obtain the following
result.

\begin{theorem}\label{th1}
The Markov chain $\{Z_n\}$  is positive recurrent if and only if inequality
\begin{equation}\label{stab2}
  r>\lambda_{1},
  \end{equation}
holds.
\end{theorem}

\noindent{\bf Remark.} Condition \eqref{stab2} may be equivalently
written as
\begin{equation}\label{stab3}
 2\gamma\lambda_1+\rho_0 -1 < \frac{(1-\rho_0)m  (\frac{\gamma^{-1}}
 {\lambda_{2}+\gamma^{-1}})^{2}
 B_{1}^{*}(\lambda_{2}) }{1+ (\frac{\gamma^{-1}}{\lambda_{2}+\gamma^{-1}} )^{2}
 B_{1}^{*}(\lambda_{2})\left[ 1- (B_{1}^{*}(\lambda_{2}))^m\right]},
  \end{equation}
where $\rho_0=\lambda_1b_1+\lambda_2b_2$. Clearly, \eqref{stab3} implies
$\rho_0<1$ and, therefore, implies \eqref{stab1}.

{\sc Proof} is based on the standard Lyapunov arguments.
We provide just a short outline of the proof.

{\it Sufficiency.}  Assume
\eqref{stab2} holds.
 Define $\tau_k$, $k=1,2,\ldots$, by
 \begin{eqnarray*}
 \tau_k&=& \inf\{n>\tau_{k-1}: Q_{n,2}=0\}, ~~ k=1,2,\ldots,
 \end{eqnarray*}
 with $\tau_0=0$. Then
$$
\Delta_k := {\mathbf E} \left((Q_{\tau_2,1} -Q_{\tau_1,1}) \ | \
Q_{\tau_1,1}=k \right) \rightarrow \lambda_1 {\mathbf E} T -
({\mathbf E} \nu +m) <0,
$$
as $k\to\infty$.  Since $\sup_{k\ge 0} \Delta_k$ is finite, the
embedded Markov chain
$Z_{\tau_k}=(Q_{\tau_k,1},Q_{\tau_k,2})=(Q_{\tau_k,1},0)$ is
positive recurrent, by the Foster criterion. Then the Markov chain
$\{Z_n\}$ is also positive recurrent, since the mean cycle time,
${\mathbf E}\tau_2-{\mathbf E}\tau_1$, is bounded from
above by ${\mathbf
E} T < \infty$.

{\it Necessity.} If \eqref{stab1} fails, then the mean time to empty the second queue
cannot be finite. Assume now that \eqref{stab1} holds and \eqref{stab2} fails.
 Then the queue lengths at the first station at the embedded moments of server's
 arrivals there couple, with a positive probability, with consecutive values of a random walk with non-negative
 drift. So the Markov chain
$\{Z_n\}$ cannot be positive recurrent. \qed

  We conclude with the following\\
  {\bf Remark.} From the expression for ${\mathbf E} \nu$, one can
  easily see that the stability condition \eqref{stab2} depends
  on the ``negative'' exponential moment $B_{1}^{*}(\lambda_{2})$,
  unless either\\
  (a) $m=0$, i.e. there is no an adaptive mechanism;
  then $r=(1-\rho )/(2\gamma + b_{1})$ or \\
  (b) $\gamma = 0$, i.e. travel times are zeros; then
  $r=(1-\rho )/b_{1}$.\\
  The stability conditions in cases (a)-(b) are well-known.

\noindent{\bf Acknowledgement.} SF thanks Onno Boxma for a useful discussion.

The research of N.\,Chernova was partially supported by the
Ministry of Higher Education and Science of the Russian Federation
Grant~RNP.2.1.1.346. The research of S.\,Foss was partially
supported by the London Mathematical Society travel grant and by
the Royal Society International Joint Project. The research of B.
Kim was supported by Basic Science Research Program through the
National Research Foundation of Korea(NRF) funded by the Ministry
of Education, Science and Technology(2011-0004133).

 \end{document}